\documentclass[12pt]{article}
\usepackage{amsmath}
\usepackage{amssymb}
\usepackage{amsthm}
\usepackage[totalwidth=17cm, totalheight=25cm]{geometry}
\usepackage{graphicx}
\usepackage{mathabx}

\newtheorem{theorem}{Theorem}[section]
\newtheorem{definition}{Definition}[section]
\newtheorem{remark}[theorem]{Remark}

\newtheorem{lemma}[section]{Lemma}

\newtheorem{thm}[equation]{Theorem}

\DeclareMathOperator{\arcsinh}{arcsinh}

\DeclareMathOperator{\Isom}{Isom}

\def\<{\langle}
\def\>{\rangle}

\def\H{\mathbb H}
\def\N{\mathbb N}
\def\R{\mathbb R}

\def\D{\partial}

\def\ga{\gamma}

\def\embed{\hookrightarrow}

\def\La{\Lambda}

\def\eps{\epsilon}
\def\la{\lambda}
\def\Ga{\Gamma}
\def\Si{\Sigma}

\def\geo{\partial_{\infty}}

\title{A note on Selberg's Lemma and negatively curved Hadamard manifolds}
\author{Michael Kapovich}

\date{\today}							

\begin{document}



\maketitle

\begin{abstract}
\noindent We prove that the conclusion of Selberg's Lemma fails for discrete isometry groups of negatively curved Hadamard manifolds. 
\end{abstract}

In  this note we give a negative answer to the first question on Margulis' problem list 
\cite[pg. 27]{Margulis}: Margulis asked if the conclusion of  Selberg's Lemma holds for  finitely generated   
isometry groups of  Hadamard manifolds.

\begin{thm}\label{thm:main}
For every $\eps>0$ and $n\ge 4$ there exists an $n$-dimensional Hadamard manifold $X_{\eps}$ of sectional curvature $-1-\eps\le K_X\le -1$ and a finitely generated discrete isometry group $\Ga_{\eps}< \Isom(X_\eps)$ which has unbounded torsion. 
\end{thm}

The idea of the proof is  simple: We start with a complete hyperbolic $n$-manifold $M^n, n\ge 4$, with finitely-generated (actually, free) fundamental group and infinitely many rank one cusps $C^n_{\la_i}$: 
Such examples were constructed in \cite{KP, K}.  
We then replace all but finitely many cusps $C^n_{\la_i}$ by metrically complete 
negatively curved (with pinching constants $(1+\eps)^{-1}$) orbifolds $O^n_i$ with boundary, where $\pi_1(O^n_i)$ is cyclic of order $i$. 
The result of this ``cusp closing'' is a complete negatively curved orbifold $O_\eps$; the action of $\Ga_\eps:= \pi_1(O_\eps)$ on the universal cover 
$X_\eps$ of $O_\eps$ provides the required examples.  

\medskip 
The Riemannian metrics in Theorem \ref{thm:main} are $C^\infty$ but not real-analytic. It is unclear if 
Theorem \ref{thm:main} holds in the real-analytic category.

\medskip
Observe that the above question has positive answer for properly discontinuous group actions in dimension $3$ 
(and, hence, $2$, although the 2-dimensional case is elementary):   Given a smooth contractible  3-manifold $X$ and a faithful properly discontinuous smooth action $\Ga\times X\to X$ of a finitely-generated group $\Ga$, 
there exists an orbifold analogue of the Scott compact core $O_c$ of the orbifold 
$O=X/\Ga$; see  \cite{FM}.  In particular, $\Ga$ is isomorphic to the fundamental group of the compact orbifold $O_c$. 
According to \cite{Hempel} and  the geometrization theorem for good compact 3-dimensional orbifolds (see \cite{BLP} or \cite{KL}), 
the orbifold $O_c$ is {\em very good}, i.e.  $\Ga$ contains a torsion-free subgroup of finite index. 

\medskip 
{\bf Acknowledgements.} The  author was partly supported by the NSF grant  DMS-16-04241, 
by the Simons Foundation Fellowship, grant number 391602,  
and  by Max Plank Institute for Mathematics in Bonn.

\section{Cusps in hyperbolic manifolds} 

We will use the upper half-space model  $\H^n=\{{\mathbf x}: x_n >0\}$ of the hyperbolic $n$-space.  An isometry of $\H^n$ is 
{\em unipotent} if it is conjugate to a translation  ${\mathbf x}\mapsto {\mathbf x} + a{\mathbf e}_1$ for  $a\ge 0$.  An isometry of $\H^n$ is called {\em parabolic} if it has a unique fixed point in the closed ball compactification $\H^n \cup \geo \H^n$ of $\H^n$.  Here $\geo \H^n$ is the {\em ideal/visual boundary} of $\H^n$. 

We let $\beta_\la: \H^n\to \R$ denote the {\em Busemann function} for the point $\la\in \geo \H^n$; this function is uniquely defined up to an additive constant.   Sublevel sets of Buisemann functions are {\em horoballs} in $\H^n$.  

Throughout the paper we will be using only closed horoballs and closed metric neighborhoods. 

Let $\H^n\embed \H^{N}$ denote an isometric totally-geodesic embedding, $N\ge n$.   
This embedding is equivariant under a  canonical monomorphism $\Isom(\H^n)\to \Isom(\H^N)$: Each isometry 
$\phi$ of $\H^n$ extends to an isometry of $\H^N$ acting trivially on the normal bundle of $\H^n$ in $\H^N$.  

For every hyperbolic subspace $X'=\H^n\subset X=\H^N$ we have the orthogonal projection $p_{X',X}:  X \to X'$. Fibers of this projection are hyperbolic subspaces orthogonal to $X'$.  In the case of nested hyperbolic subspaces 
$$
X''\subset X'\subset X, 
$$
we have
\begin{equation}\label{eq:nested}
p_{X'',X}= p_{X'',X'} \circ p_{X',X}.  
\end{equation}

\medskip
Below we review the notion of {\em cusps} of hyperbolic manifolds/orbifolds; we refer to 
\cite{Bowditch93, Bowditch95, Rat} for details.  

Let $\Ga< \Isom(\H^n)$ be a discrete subgroup with the limit set $\La=\La(\Ga)\subset \D_{\infty} \H^n$.  
A {\em parabolic limit point} of $\Ga$ is a fixed point of a parabolic isometry  
$\ga\in\Ga$. The $\Ga$-stabilizer $\Pi=\Ga_\la< \Ga$ of such $\la\in \La$ is called  {\em a maximal parabolic subgroup} of $\Ga$. 
 A parabolic limit point $\la$ of $\Ga$ is called {\em bounded} (equivalently, {\em cusped}) if the quotient $(\La-\{\la\})/\Ga_\la$ 
is compact. Each bounded parabolic fixed point of $\Ga$ corresponds to a 
``cusp'' of the quotient orbifold $M=\H^n/\Ga$ defined as follows. Let $X'_\la\subset X=\H^n$ be a smallest $\Pi$-invariant hyperbolic subspace of $X$ (such a subspace need not be unique). Then $\Pi$ acts with finite covolume on the intersection $B_\la\cap X'_\la$ for every horoball $B_\la\subset \H^n$ centered at $\la$.  The virtual rank of the virtually abelian group $\Pi$ equals $r_\la=\dim(X'_\la)-1$. 

\medskip 
Let $p_\la:= p_{X'_\la,X}:  X=\H^n\to X'_\la$ be the orthogonal projection as above.   
Define 
$$
\tilde{C}'_\la:= B_\la\cap X'_\la
$$
and 
$$
\tilde{C}_\la:= p_\la^{-1}(\tilde{C}'_\la)\subset \H^n$$
 (both depend on $B_\la$ and $X'_\la$, of course). 

\begin{definition}
If the orbi-covering map $\pi: \H^n/\Pi\to \H^n/\Ga=M$ is injective on $\tilde{C}_\la/\Pi$, then the image 
${C}_\la  := \pi( \tilde{C}_\la/\Pi )$  is called a  {\em cusp neighborhood} in $M$ (or, simpy, a {\em cusp} in $M$) corresponding to $\Pi$.  
The domain $\tilde{C}_\la$ is then called a {\em cusped region} of the limit point $\la\in \La(\Ga)$.  
The number $r_\la$ is the {\em rank} of the cusp $C_\la$. 
\end{definition}

By abusing the notation, for a cusped region  $\tilde{C}_\la$,  
we will denote $C_\la=\tilde{C}_\la/\Pi $ as well. 

\medskip 
For each $n$-dimensional cusp $C_\la$ we define its  {\em core}  $C'_\la\subset C_\la$ as the quotient $\tilde{C}'_\la/\Pi$.  
The core is unique up to an isometry $C_\la\to C_\la$.

\medskip 
A parabolic limit point $\la\in \La(\Ga)$ is a bounded if and only if  for a sufficiently small horoball $B_\la$ 
(depending, among other things, on the choice of $X'_\la$)  $\tilde{C}_\la$ is a cusped region.

\begin{remark}
 Each maximal parabolic subgroup $\Ga_\la<\Ga$ (regardless of whether $\la$ is a bounded parabolic limit point or not) 
corresponds to a {\em Margulis cusp} of $M=\H^n/\Ga$: It is the projection to $M$ of the region $U_\la\subset \H^n$ consisting of 
points $x$ such that there exists a parabolic element $\ga\in\Ga_\la$ satisfying $d(x, \ga(x))\le \mu_n$, the Margulis constant of $\H^n$. 
  Margulis cusps should not be confused with the cusps defined above. Margulis cusps will not be used in this paper.   
\end{remark}

If  $\Ga< \Isom(\H^n)$ is {\em geometrically finite} then every parabolic limit point of $\Ga$ is bounded, $M$ has only finitely many cusps 
and, after taking sufficiently small horoballs $B_\la$, we can assume that these cusps are pairwise  disjoint.    If $n=3$ then every finitely-generated discrete subgroup $\Ga< \Isom(\H^3)$ has only finitely many  cusps (and Margulis cusps), 
but this fails in dimensions $n\ge 4$  (see \cite{KP});  the existence of such examples is critical for the proof of Theorem \ref{thm:main}. 

\medskip
We will call  a cusp {\em unipotent} if its fundamental group $\Pi$ 
is unipotent, i.e.  every element of $\Pi$ is unipotent.

\medskip
For each $r>0$ we define the {\em $r$-collar} $C_{\la,r}\subset C_\la$ of the boundary of a cusp $C_\la$ as the quotient by $\Pi$ of $p_\la^{-1}(N_r( X'_\la \cap \D  B_\la))$, where $N_r(\cdot)$ denotes the  $r$-neighborhood in $X'_\la \cap B_\la$.  
Then the minimal distance between the boundary components of the collar $C_{\la,r}$ equals $r$. 

\medskip 
Given a hyperbolic subspace $\H^n\subset \H^{N}$, every horoball $B_\la\subset \H^n$ properly embeds in a horoball $B^N_\la\subset \H^N$.  Accordingly, each cusp $C^n_\la\subset \H^n/\Pi$ properly embeds in a cusp $C^N_\la\subset \H^N/\Pi$.

\medskip
In this paper we will be only interested in rank one unipotent cusps $C_\la$. Such a cusp is uniquely determined  
(up to isometry)  by its dimension $n$ and one more real parameter, the {\em core length} $\ell(C_\la)$, defined as the length of the 
boundary loop of the 2-dimensional core $C'_\la\subset C_\la$. 

\medskip 
We will need the following example of a finitely generated discrete subgroup of $\Isom(\H^4)$ with infinitely many cusps:

\begin{thm}
[\cite{KP, K}] 
\label{thm:KP}
There exists a discrete (geometrically infinite) subgroup $\Phi< \Isom(\H^4)$ isomorphic to a free group $F_k$ of rank $k<\infty$, such that:

1.  The quotient manifold $M^4=\H^4/\Phi$ contains an infinite collection of  pairwise disjoint and  isometric  
rank 1  unipotent cusps $C_{\la_i}, i\in \N$.

2. $\Phi$ is a normal subgroup of a geometrically finite group $\hat\Phi< \Isom(\H^4)$; every 
rank one cusp of $M^4$ injectively covers a rank one  cusp of $\H^4/\hat\Phi$.  
\end{thm}

We let $L:= \ell(C_{\la_i})$ denote the common core length of the cusps $C_{\la_i}$ of  $M$.  

\begin{remark}
Besides the cusps $C_{\la_i}$, the manifold $M^4$ also has finitely many Margulis cusps. 
These Margulis cusps project  to rank two cusps of $\H^4/\hat\Phi$. 
The parabolic limit points of $\Phi$ corresponding to these Margulis cusps are not bounded. \end{remark}

\medskip 

We retain the notation $\Phi$ for the image of $\Phi$ under the embedding 
$\Isom(\H^4)\to \Isom(\H^n)$ and define $M^n:= \H^n/\Phi$.  As noted above, each cusp $C^4_{\la_i}$ of 
$M^4$ embeds properly in a cusp $C_{\la_i}^n\subset M^n$; these $n$-dimensional cusps are also pairwise disjoint, since 
there exist $1$-Lipschitz retracts $M^n\to M^4$, satisfying $C^n_{\la_i}\to C^4_{\la_i}$.

\section{Warped products} 

The ``cusp closing'' procedure in the proof of Theorem \ref{thm:main} will use {\em warped products} of negatively curved manifolds. In this section we review this basic construction.

Let $(B, ds_B^2), (F, ds^2_F)$ be Riemannian manifolds and $f: B\to (0,\infty)$ a smooth function. The {\em warped product} of these manifolds, denoted $W=B\times_f F$, is the product $B\times F$ equipped with the Riemannian metric 
$$
ds^2= ds_B^2 + f^2 ds^2_F, 
$$
see  \cite[sect. 7]{BO}  for a detailed discussion. In particular, the Riemannian manifold $W$ is complete iff $B$ and $F$ 
both are; \cite[Lemma 7.2]{BO}. 

If $\Gamma$ is a group acting isometrically on each factor $(B, ds_B^2), (F, ds^2_F)$  and $f: B\to \R$ is $\Ga$-invariant, 
then the product action of $\Gamma$ on $W$ is also isometric.  In particular, the notion of warped product extends to 
Riemannian orbifolds. 

We will need the sectional curvature formula for the warped product  in the case of 
$(B, d_B^2)=(\R, dt^2)$, given in \cite[pg. 27]{BO}:   
$$
K(\Pi)= -\frac{f''(t)}{f(t)} ||x||^2 + \frac{L(v,w) - (f'(t))^2}{f^2(t)}||v||^2.
$$
Here $\Pi$ is  plane in $T_{(t,q)}W$ with the orthonormal basis $\{x+v, w\}$, where $v, w\in T_qF$, $x$
 is a horizontal tangent vector (thus, $||x||^2 + ||v||^2= ||x||^2 + f^2 ||v||_F^2=1$)   
and $L(v,w)$ is the sectional curvature of $(F, ds^2_F)$ at $q$ on the plane spanned by the vectors $v, w$.  (Note that  \cite[pg. 26]{BO} 
also contains the sectional curvature formula for general warped products).

In particular, if $F$ is negatively curved with sectional curvature $-1-\eps \le L\le -1$ and $f(t)=\cosh(t)$, then 
$$
-1-\eps\le  K(\Pi)\le -1
$$
as well. 

The hyperbolic space $\H^n$ is isometric to a warped product $\H^{n-2} \times_f \H^{2}$ with $f: \H^{n-2}\to \R_+$, $f(p)= \cosh(d(o, p))$, where $o\in \H^{n-2}$ is a basepoint. This warped product decomposition can be realized as follows. We let $\H^{2}$ be embedded in $\H^n$ as
$$
\{(x_1, 0, 0, ....,0, x_n): x_n > 0\}. 
$$
Horizontal (totally-geodesic) leaves of the warped product $\H^{n-2} \times_f \H^{2}$ correspond to 
codimension two hyperbolic subspaces in $\H^n$ orthogonal to $\H^{2}$, 
while vertical leaves are obtained from  $\H^{2}$ by rotating it via elements of $SO(n-1)< SO(n)$ fixing pointwise the coordinate line 
$$
\R {\mathbf e}_n = \{(x_1, 0,  ...., 0, 0)\}.
$$ 
The vertical projection $\H^{n-2} \times_f \H^{2}\to \H^2$ is just the orthogonal projection 
$p_{\H^2,\H^n}$.

Yes another way to realize this decomposition of $\H^n$ is as the iterated warped product 
$$
\H^n= \R\times^{n-2}_f \H^2,
$$ 
$$
\H^3= \R\times_f \H^2, \H^4= \R\times_f \H^3, ..., \H^n= \R\times_f \H^{n-1}, 
$$
where $f(t)=\cosh(t)$.  The orthogonal  projection $p_{\H^2,\H^n}$ equals the vertical projection $\eta: \R\times^{n-2}_f \H^2\to \H^2$ 
given by iterating vertical projections in the warped product decompositions $\H^k=\R \times_f \H^{k-1}$, cf. \eqref{eq:nested}.

\medskip 
We generalize this iterated warped product as follows. We let $\tilde{F}$ be a simply-connected complete 
negatively  curved surface with sectional curvature in the interval $[-1-\eps, -1]$. Define the iterated warped product 
\begin{equation}\label{eq:iwp}
\tilde{W}=    \R\times^{n-2}_{\cosh} \tilde{F}. 
\end{equation}
 It follows that $\tilde{W}$ is still  a simply-connected complete 
negatively  curved manifold with sectional curvature in $[-1-\eps, -1]$.  We let $\tilde\eta: \tilde{W}\to \tilde{F}$ 
denote the vertical projection. Then for an open subset $\tilde{U}\subset \tilde{F}$, the preimage $\tilde\eta^{-1}(\tilde{U})$ 
is  an iterated warped product $\R  \times^{n-2}_f \tilde{U}$. In particular, if $\tilde{U}$ 
has constant curvature $-1$, so does $\tilde\eta^{-1}(\tilde{U})$.

\section{Closing rank one cusps}

We will apply  iterated warped products \eqref{eq:iwp} to surfaces $\tilde{F}$ (and their Riemannian orbifold quotients), 
constructed by splicing quotients of $\H^2$ by cyclic parabolic and by finite cyclic groups.  
The goal is  ``close'' $n$-dimensional rank one unipotent cusps $C^n_\la$, converting them to orbifolds of variable negative curvature  
with finite cyclic fundamental groups, while leaving the Riemannian metric on a suitable $r$-collar of $C_\la$ unchanged.   The cusp-closing 
is a rather standard procedure, we describe it here in detail for the sake of completeness. 

We start by describing  cusp-closing in dimension 2.  
Let $\Si_0< \Isom(\H^2)$ be a cyclic parabolic subgroup; the surface $T_0:= \H^2/\Si_0$ is foliated by projections of $\Si_0$-invariant horocycles in $\H^2$. Let $c_0\subset T_0$ be the (unique) leaf of length $a >0$.  (The number $a$ will be specified later on.)

Similarly, let $\Sigma_i< \Isom(\H^2)$ be a  finite  cyclic subgroup of order $i\ge 2$. The quotient-orbifold $T_i=\H^2/\Sigma_i$ is foliated by projections of $\Sigma_i$-invariant circles in $\H^2$. 
Let $c_i\subset T_i$ be the (unique)  
leaf of  the same length $a$ as above.   The hyperbolic  
surfaces/orbifolds $T_0$ and $T_i$ admit isometric $U(1)$-actions whose orbits are 
leaves of the above foliations.  The distance from $c_i$ to the singular point of $T_i$ 
equals $R_i=\arcsinh(\frac{ai}{2\pi})$. 

We let $T'_0$ denote the closure of the infinite 
 area component in $T_0-c_0$ and let $T''_i$ denote the closure of the bounded component of 
$T_i- c_i$. Gluing $T'_0, T''_i$ via an isometry of their boundaries results in a metric orbifold $S_i$; the metric on $S_i$ is, of course, smooth away from $\bar{c}_i:= c_0\equiv c_i$ and is singular along that curve.  (The group $U(1)$ still acts isometrically on $S_i$.) 
Below we  smooth out the metric on $S_i$ by modifying it near 
$\bar{c}_i$, so that the new metric has negative curvature with small pinching constant when $i$ is large.

Fix $r>0$ and let the nested annuli $A_{i,r}\subset A_{i,2r}$ denote the $r$- and $2r$-neighborhoods of $\bar{c}_i$ in $S_i$ with respect to the singular  metric on $S_i$.  We will take $r$ such that 
$$
r< \frac{1}{2} \arcsinh(a/\pi)\le R_i,
$$
hence, the annulus $A_{i,2r}$ is disjoint from the singular point of the orbifold $S_i$.  

The next lemma follows from the geometric convergence of suitable conjugates of subgroups $\Si_i< \Isom(\H^2)$ to 
$\Si_0< \Isom(\H^2)$, since the latter implies $C^\infty$ Gromov--Hausdorff  convergence 
$$
(T_i, t_i)\to (T_0,t_0),$$ 
where $t_i\in c_i, t_0\in c_0$; see \cite[Ch. E]{BP}.

\begin{lemma}\label{lem:L}
For each $\eps> 0$ there is $i_\eps$ such that for all $i\ge i_\eps$ 
there exist  $U(1)$-invariant Riemannian metrics $g_i$ on the orbifolds $S_i$ satisfying:

1. $g_i$ equals the restrictions of the metrics of $T_0$ and $T_i$ respectively on the unbounded/bounded components of $S_i -A_{i,r}$. 

2. The curvature of $g_i$ lies in the interval $[-1-\eps, -1]$.    
\end{lemma}

In what follows, we equip the orbifolds $S_i$ with the above metrics $g_i$ and denote the resulting Riemannian orbifold $F_i$.  We let $\tilde{F}_i\to F_i$ denote the (degree $i$) universal cover of the orbifold $F_i$; then $\tilde{F}_i$ is a  simply-connected 
negatively curved complete Riemannian surface. We let $O^2_{i,2r}\subset F_i$ denote the union of the annulus 
$A_{i,2r}\subset F_i$ and the suborbifold $T''_i\subset F_i$ (equipped with the restriction of the metric $g_i$, of course).  The boundary curve of $O^2_{i,2r}$ has length 
$a e^{2r}$. 

 \medskip 
Before extending this construction to higher dimensions we describe  cusp-closing for hyperbolic surfaces. Let $M^2$ be a complete hyperbolic surface (possibly of infinite area) and $C_{\la_i}\subset M^2$ be pairwise disjoint cusps with equal core lengths $= L > a$.  Each $C_{\la_i}$, of course embeds isometrically in the surface $T_0$ as above; we let $2r$ denote the distance between the boundary curve of $C_{\la_i}$ and the loop $c_0\subset T_0$.  Thus, $L=a e^{2r}$.

The $r$-collar $C_{\la_i,r}\subset C_{\la_i}$ is isometric to the $r$-neighborhood 
$E_{i,r}=N_r(\D O^2_{i,2r})$ of the boundary 
of $O^2_{i,2r}$ ($E_{i,r}$ is a component of  $A_{i,2r} - A_{i,r}$ and has constant negative curvature).  
Hence, we can replace each cusp $C_{\la_i}\subset M^2$ 
with a Riemannian orbifold $O^2_{i,2r}$ by first removing $C_{\la_i} - C_{\la_i,r}$ and then gluing $O^2_{i,2r}$ via an isometry 
$E_{i,r}\to C_{\la_i,r}$.  The resuling Riemannian orbifold $O^2$ is said to be obtained from $M^2$ by {\em cusp-closing}.  

\begin{remark}
The Riemannian orbifold $O^2$ is complete and has negative curvature, which equals $-1$ except for the annuli $A_{i,r}$ where the curvature is variable. In dimension $2$ (at least if the group $\pi_1(M^2)$ is finitely generated) one can  also accomplish cusp-closing via a metric of constant negative curvature (by perturbing the hyperbolic metric of $M^2$ globally rather than inside of the cusps $C_{\la_i}$), but this is not what we are interested in.  See also the \cite[Corollary 2]{K} in the setting of hyperbolic $4$-manifolds/orbifolds. 
\end{remark}

We now proceed with cusp-closing in higher dimensions. 

\medskip 
Applying the iterated warped product construction to the Riemannian surfaces $\tilde{F}= \tilde{F}_i$ 
as described above, we obtain $n$-dimensional Hadamard manifolds $\tilde{W}^n_i$ equipped with isometric $\Si_i$-actions; let 
$W^n_i:= \tilde{W}^n_i/\Si_i$ be the Riemannian quotient-orbifolds. 
 Thus, each $W^n_i$ is the $n-2$-fold iterated warped product
$$
W^n_i= \R \times_{\cosh}^{n-2} F_i
$$
with the vertical projection $\eta_i: W_i\to F_i$.  As noted above, $W_i$ has constant curvature $-1$ away from $\eta^{-1}(A_{i,r})$. 

As with the 2-dimensional cusp closing, we will only need the parts  
$$
O^n_{i,2r} := \eta^{-1}(O^2_{i,2r})\subset W^n_i$$
of the orbifolds $W^n_i$.  Each $O^n_{i,2r}$ can be regarded as an $\H^{n-2}$-bundle over the orbifold $O^2_{i,2r}$. 
 The orbifolds  $O^n_{i,2r}$  will be replacing rank one unipotent cusps of a hyperbolic $n$-manifold. This will be accomplished by 
gluing along constant curvature boundary collars 
$$
E^n_{i,r}:= \eta^{-1}(E_{i,r})\subset O^n_{i,2r}. 
$$

Let $C^n_{\la}$ be an $n$-dimensional rank one unipotent cusp of core length $\ell(C^n_\la)=L=a e^{2r}$.  In view of the iterated warped product decomposition 
$$
\H^n= \R \times^{n-2}_{\cosh} \H^2, 
$$
the cusp $C^n_{\la}$ also decomposes as the iterated warped product
$$
C^n_{\la}= \R \times^{n-2}_{\cosh} C^2_\la,
$$
where $C^2_\la\subset C^n_\la$ is a 2-dimensional core of $C^n_\la$, 
with the projection $\eta: C^n_{\la}\to C^2_\la$.   The boundary $r$-collar $C^n_{\la,r}\subset C^n_{\la}$ equals the preimage 
$$
\eta^{-1}(C^2_{\la,r})
$$
of the $r$-collar of the core cusp $C^2_\la$. 

Thus we obtain isometries of boundary collars 
$$
C^n_{\la,r} \to E^n_{i,r},
$$
from the boundary collar of a cusp $C^n_\la$ to the boundary collars of the orbifolds $ O^n_{i,2r}$.

\section{Proof of Theorem \ref{thm:main}}

We construct a Riemannian orbifold $O^n_\eps$ as follows. Given $\eps> 0$, we let $i_\eps\in \N$ be as in Lemma \ref{lem:L}.  
Recall that $L=\ell(C^4_{\la_i})$ is the common core length of the cusps $C^4_{\la_i}$ of the hyperbolic $4$-manifold $M^4$ 
from Theorem \ref{thm:KP}.  We then let $r>0$ and $a>0$ be such that 
$$
L= e^{2r} a, \quad r< \frac{1}{2} \arcsinh(a/\pi), 
$$ 
which can be always accomplished by taking $r$ to be sufficiently small.

 From each cusp $C^n_{\la_i}, i\ge i_\eps$, of the hyperbolic $n$-manifold $M^n=\H^n/\Phi$ we remove the complement to the boundary 
$r$-collar $C^n_{\la_i,r}$; let $M'$ denote the remaining manifold:
$$
M':= M^n - \coprod_{i\ge i_\eps} (C^n_{\la_i} -  C^n_{\la_i,r}). 
$$
Then for each $i\ge i_\eps$ we glue to  $M'$ the Riemannian orbifold $O^n_{i,2r}$ via an isometry of the collars
$$
C^n_{\la_i,r} \to E^n_{i,r},
$$
The result is an $n$-dimensional  Riemannian orbifold $O_\eps:= O^n_\eps$.  

\begin{remark}
Since $\pi_1(M^n)\cong \pi_1(M^4)\cong \Phi$ is free of rank $k$, the fundamental group $\Ga_\eps=\pi_1(O_\eps)$ has the presentation
$$
\< s_1,...,s_k| w_i^i, i\ge i_\eps\>,
$$
where the words $w_i$ represent  generators of fundamental groups of the cusps $C^4_{\la_i}\subset M^4$. 
\end{remark}

By the construction, the sectional curvature of $O_\eps$ lies in the interval $[-1-\eps, -1]$.  Since  $M'$ and all orbifolds   
$O^n_{i,2r}$ are metrically complete and the minimal distance between the boundary components of each collar 
$C^n_{\la_i,r}$ equals $r>0$, it follows that the  Riemannian orbifold $O_\eps$ is also complete.

Since the orbifold $O_\eps$ is complete and negatively curved, it is good (developable); see 
\cite[pg. 603, Theorem 2.15]{BH} and also \cite[Ch. 13]{Rat}.   
Hence, the universal cover of $O_\eps$ is an $n$-dimensional 
Hadamard manifold $X=X_\eps$  of curvature $-1-\eps_m \le K_X\le -1$. The fundamental 
group $\Ga_\eps$ acts on $X_\eps$ faithfully, properly discontinuously and isometrically with 
$O_\eps \cong X_\eps/\Ga_\eps$. In particular, $\Ga_\eps$ has unbounded torsion: 
The fundamental group (cyclic group of order $i$) of each 
orbifold $O^n_{i,2r}$ ($i\ge i_\eps$) embeds in $\Ga_\eps$.   \qed

\noindent Department of Mathematics, UC Davis, One Shields Avenue, Davis, CA 95616, USA\\
kapovich@math.ucdavis.edu

\end{document}